\documentclass[draft]{amsart}
\usepackage[active]{srcltx}
\usepackage{atveryend}
\makeatletter
\let\origcitation\citation
\AtEndDocument{\def\mycites{}
  \def\citation#1{\g@addto@macro\mycites{#1^^J}\origcitation{#1}}}
\AtVeryEndDocument{\newwrite\citeout\immediate\openout\citeout=\jobname.cit
  \immediate\write\citeout{\mycites}\immediate\closeout\citeout}
\makeatother
\usepackage{enumerate}
\usepackage{cancel}
\usepackage{mathdots}
\usepackage{epsfig}
\usepackage{hyperref}
\usepackage{amsmath}
\usepackage{amssymb}
\usepackage{graphics}
\usepackage[Symbol]{upgreek} 

\usepackage{yfonts}
\newtheorem*{Theorem*}{Theorem}

\newtheorem*{Lemma 1*}{Lemma 1}
\newtheorem*{Lemma 2*}{Lemma 2}
\newtheorem*{Lemma 3*}{Lemma 3}
\newtheorem*{Lemma 4*}{Lemma 4}

\DeclareFontFamily{U}{fsy}{}
\DeclareFontShape{U}{fsy}{m}{n}{<->s*[.9]psyr}{}
\DeclareSymbolFont{der@m}{U}{fsy}{m}{n}
\DeclareMathSymbol{\der}{\mathord}{der@m}{182}

\def\dbm {\der_{\operatorname{BM}}}
\def\No {\mathbf{No}}
\def\cf {\operatorname{cf}}

\title{homogeneous universal \emph{H}-fields}


\author{Lou van den Dries}

\address{Department of Mathematics\\
University of Illinois at Urbana-Champaign\\
Urbana, IL 61801, USA}

\email{vddries@illinois.edu}

\author{Philip Ehrlich}

\address{Department of Philosophy\\
Ohio University\\
Athens, OH 45701, USA}

\email{ehrlich@ohio.edu}

\begin{document}

\begin{abstract} 
We consider derivations $\der$ on Conway's field $\mathbf{No}$ of surreal numbers such that the ordered differential field $(\mathbf{No},\der)$ has constant field $\mathbb{R}$ and is a model of the model companion of the theory of $H$-fields with small derivation.
We show that this determines $(\mathbf{No},\der)$ uniquely up to isomorphism, and that this structure is
absolutely homogeneous universal for models of this theory with constant field $\mathbb{R}$.

\end{abstract}

\maketitle

\noindent
Aschenbrenner and van den Dries~\cite{AD} introduced $H$-fields in order to formalize some basic first-order properties of Hardy fields in their role of ordered and valued differential fields. Hardy fields containing $\mathbb{R}$ are $H$-fields, and so is the system $\mathbb{T}$ of \emph{transseries}. In \cite{ADH1}, Aschenbrenner, van den Dries, and van der Hoeven (ADH) proved that the theory of $H$-fields has a model companion whose models are the $H$-fields
that are Liouville closed, $\upomega$-free, and newtonian.
Adding to these axioms for the model companion the axiom that the derivation is small yields a complete theory $T$ in the language $\mathcal {L}=\left \{ 0,\ 1,\ +,\ \cdot,\ \le,\ \preceq,\ \der \right \}$ of ordered valued differential fields. Thus $T$ is complete as well as model complete. Another key result from \cite{ADH1} is that $\mathbb{T}$ is a model of $T$. See
\cite[Introduction]{ADH1} for the relevant definitions.

{ }\quad Using an idea from Schmeling's thesis~\cite{S} due to van der Hoeven, Berarducci and Mantova \cite{B} constructed so-called \emph{surreal derivations} on Conway's ordered field $\mathbf{No}$ of \emph{surreal numbers} \cite{C}; even  a ``simplest" one, $\dbm$, that makes 
$\big(\mathbf{No}, \dbm \big)$
an $H$-field with small derivation and constant field $\mathbb{R}$. They proved also that this $H$-field is Liouville closed. ADH \cite{ADH2} subsequently showed that $\big(\mathbf{No}, \dbm \big)$ is a model of $T$ that is universal with respect to $H$-fields with small derivation and constant field $\mathbb{R}$: every such $H$-field, including each Hardy field containing $\mathbb{R}$, can be embedded as an ordered differential field into 
$\big(\mathbf{No}, \dbm\big)$.  The purpose of this note is to point out that in the course of establishing the just-said result, \cite{ADH2} proves almost enough to obtain the following:

\begin{Theorem*}\label{1} Let $\der$ be any derivation on $\mathbf{No}$
with constant field $\mathbb{R}$ such that $\big(\mathbf{No}, \der \big)$ is a model of $T$. Then $\big(\mathbf{No}, \der \big)$
is up to isomorphism the unique model of $T$ with constant field 
$\mathbb{R}$ that is absolutely homogeneous universal with respect to models of $T$ with constant field $\mathbb{R}$.  
\end{Theorem*}
 
\noindent
The uniqueness gives $\big(\mathbf{No}, \der \big)\cong
\big(\mathbf{No}, \dbm \big)$. Part of the interest of the theorem lies
in the circumstance that $\dbm$ seems to take the ``wrong'' values on some infinite iterates of the exponential function applied to $\omega$; for more on this, see \cite{ADH3, B1}.  We do expect there is an {\em optimal\/} surreal derivation--- better than $\dbm$---that also satisfies the hypothesis of the theorem, and thus its conclusion. 

Let $\der$ be as in the theorem. That $\big(\mathbf{No}, \der \big)$ is \emph{absolutely universal with respect to models of $T$ with constant field $\mathbb{R}$\/} means that every model of $T$ with constant field $\mathbb{R}$ (whose universe is a set or a proper class) can be embedded in 
$\big(\mathbf{No}, \der \big)$. That it is \emph{absolutely homogeneous with respect to models of $T$ with constant field $\mathbb{R}$\/} means that every isomorphism between substructures of $\big(\mathbf{No}, \der \big)$ that are set-models of $T$ with constant field $\mathbb{R}$ extends to an automorphism of $\big(\mathbf{No}, \der \big)$. 

Our set theory here is von Neumann-Bernays-G\"odel set theory with Global Choice (NBG), a conservative extension of ZFC in which all proper classes are in bijective correspondence with the class $On$ of all ordinals. By ``set-model" (``class-model") we mean a model whose universe is a set (a proper class). By ``cardinal'' we mean below ``set-cardinal'' and we let $\kappa$ range over cardinals.

\medskip\noindent
To establish the theorem we follow the proof of \cite[Theorem 3]{ADH2},
which deals with the case $\der=\dbm$. To handle arbitrary $\der$ we use some extra lemmas. The first one slightly extends \cite[Lemma 5.3]{ADH2}, which considers only regular $\kappa$. The proof of that lemma goes through if we replace $\kappa$ at various places by its cofinality
$\cf(\kappa)$.

\begin{Lemma 1*}\label{lemsat} Let $\kappa$ be uncountable. Then the underlying ordered sets of $\mathbf{No}(\kappa)$ and 
$v\big(\No(\kappa)^\times\big)$ are $\cf(\kappa)$-saturated.  
\end{Lemma 1*}

\begin{Lemma 2*}\label{cfL} Let $L$ be a countable $($one-sorted$)$ language and $\mathbf{No}_L$ an $L$-structure with universe $\mathbf{No}$. Then there are cardinals $\kappa$ of arbitrarily large cofinality such that $\mathbf{No}(\kappa)$ is the underlying set of an elementary substructure of  $\mathbf{No}_L$.
\end{Lemma 2*}
\begin{proof} By Skolemizing we arrange that $\text{Th}(\mathbf{No}_L)$ has built-in Skolem functions, so any substructure of $\mathbf{No}_L$ is an elementary substructure. Let $\kappa$ be an infinite regular cardinal. 
We build in the usual way simultaneously by transfinite recursion a strictly increasing sequence 
$(\kappa_{\alpha})_{\alpha<\kappa}$ of infinite cardinals and an elementary chain $(K_{\alpha})_{\alpha<\kappa}$
of elementary substructures of $\mathbf{No}_L$ such that for all $\alpha< \kappa$, \begin{enumerate}
\item  
$\mathbf{No}(\kappa_\alpha)\subseteq K_\alpha\subseteq \mathbf{No}(\kappa_{\alpha+1})$, where $K_{\alpha}$ denotes also its underlying set.
\item if $\alpha$ is an infinite limit ordinal, then 
$\kappa_{\alpha}=\sup_{\beta<\alpha} \kappa_{\beta}$ and
$K_{\alpha}=\bigcup_{\beta< \alpha} K_{\beta}$.
\end{enumerate}
Then $\mathbf{No}(\kappa_{\infty})$ with $\kappa_{\infty}:= \sup_{\alpha< \kappa} \kappa_\alpha$
is the underlying set of the elementary substructure 
$\bigcup_{\alpha<\kappa} K_{\alpha}$ of $\mathbf{No}_L$, and $\kappa_{\infty}$
has cofinality $\kappa$.  
\end{proof} 

\noindent
In the next two lemmas ``$H$-field'' should be read as ``$H$-field whose universe is a set''. We note that any embedding between $H$-fields with common constant field $\mathbb{R}$ is automatically the identity on $\mathbb{R}$. 
In the rest of the paper we fix a class-model $\big(\mathbf{No}, \der \big)$ of $T$ with constant field $\mathbb{R}$. Here is the relevant analogue of the Claim in the proof of \cite[Theorem 3]{ADH2}:

\begin{Lemma 3*}\label{claim}  Let $E\subseteq K$ be an extension of $\upomega$-free $H$-fields with $\mathbb{R}$ as their common constant field, and let $i: E \to (\mathbf{No},\der)$ be an embedding. Then $i$ extends to an embedding $K \to (\mathbf{No},\der)$. 
\end{Lemma 3*} 
\begin{proof} First extend $K$ to make it a model of $T$; by \cite[16.4.1 and 14.5.10]{ADH1} this can be done 
without changing its constant field. By Lemma 2 we can take an
uncountable cardinal 
$\kappa$ such that $\cf(\kappa)> \text{card}(K)$, $\mathbf{No}(\kappa)$
underlies an elementary substructure $L$ of $(\mathbf{No},\der)$
and $i(E)\subseteq \mathbf{No}(\kappa)$. Using 
Lemma 1 and \cite[16.2.3]{ADH1} we then extend $i$ 
to an embedding $K \to L$.
\end{proof}

\begin{Lemma 4*}\label{initial} There is an $\upomega$-free $H$-field with
constant field $\mathbb{R}$ that embeds into every model of $T$ with constant field $\mathbb{R}$.
\end{Lemma 4*}
\begin{proof} Let $F$ be the Hardy field $\mathbb{R}(x)$ (so $x>\mathbb{R}$, $x'=1$). Then
$F$ is a grounded $H$-field with constant field $\mathbb{R}$ that embeds into every
model of $T$ with constant field $\mathbb{R}$. Hence $F_{\upomega}$ as in
\cite[Lemma 11.7.17]{ADH1} is an $\upomega$-free $H$-field with
constant field $\mathbb{R}$ that embeds into every model of $T$ with constant field $\mathbb{R}$.
\end{proof}

\begin{proof}[Proof of the Theorem] Recall that $(\mathbf{No},\der)$ is a model of $T$ with constant field $\mathbb{R}$.

As to universality for set-models,
let $K$ be a set-model of $T$ with constant field $\mathbb{R}$. Use Lemma 4 to make
$K$ an extension of an $\upomega$-free $H$-field $E$ with an embedding
$E \to (\mathbf{No},\der)$, and then use Lemma 3 to extend this embedding to an embedding $K\to (\mathbf{No},\der)$. As to universality for class-models, let $K$ be a class-model of $T$ with constant field $\mathbb{R}$.
Then $K$ is the union of a chain $(K_\beta)_{\beta \in  On}$ of set-models of $T$ with constant field $\mathbb{R}$. First embed $K_0$ into $(\mathbf{No},\der)$, and then use transfinite recursion, Lemma 3, and Global Choice to construct a family $(i_\beta)_{\beta\in On}$ of embeddings $i_\beta: K_\beta\to (\mathbf{No},\der)$,
with $i_{\beta}$ extending $i_\alpha$ whenever $\alpha < \beta$. Then the common extension of these $i_\beta$ is an embedding $K\to (\mathbf{No},\der)$.

For homogeneity for set-models, 
  let $i: E \to F$ be an isomorphism between set-models
$E, F\preceq (\mathbf{No},\der)$. Given any $a\in \mathbf{No}\setminus E$
we use Lemma 3 to extend $i$ to an isomorphism $K\to L$ between set-models $K, L\preceq (\mathbf{No},\der)$ with $a\in K$. Likewise for $b\in \mathbf{No}\setminus F$
we can extend $i$ to an isomorphism $K\to L$ between set-models $K, L\preceq (\mathbf{No},\der)$ with $b\in L$. The usual back-and-forth argument then gives an automorphism of $(\mathbf{No},\der)$ extending $i$.  We have now shown that $(\mathbf{No},\der)$ is absolutely homogeneous universal with respect to models of $T$ with constant field $\mathbb{R}$.

As to uniqueness, let $M$ be any class-model of $T$ with constant field $\mathbb{R}$ that is absolutely homogeneous universal with respect to models of $T$ with constant field 
$\mathbb{R}$;  we need to show that $M\cong
\big(\mathbf{No}, \der \big)$. First, let a set-model $E\preceq M$ and an isomorphism $i: E \to F\preceq \big(\mathbf{No}, \der \big)$ be given.  
To go {\em forth},  Lemma 3 allows us to extend $i$ for any $a\in M\setminus E$ 
to an isomorphism $K\to L$ between set-models $K\preceq M$ and $L\preceq (\mathbf{No},\der)$ with $a\in K$.  To go  {\em back}, let $b\in \mathbf{No}\setminus F$,
take a set-model $L\preceq \big(\mathbf{No}, \der \big)$ with $F\subseteq L$ and $b\in L$, and take an embedding $j: L \to M$. Then $j\circ i$ maps $E$ isomorphically onto $j(i(E))\preceq M$, and so extends to an automorphism $\sigma$ of  $M$. Then $\sigma^{-1}\circ j$ extends $i^{-1}$ and maps $L$ isomorphically onto some $K\preceq M$ with $K\supseteq E$. Thus back-and-forth yields an isomorphism  $M \to \big(\mathbf{No}, \der \big)$. 
\end{proof} 

\noindent
Alternatively (but there is really little difference with the approach above) we could have adapted familiar arguments of J\'onsson \cite{J1, J2} for the existence and uniqueness (up to isomorphism) of an $L$-structure of inaccessible power $\kappa$ that is $\kappa$-homogeneous and $\kappa$-universal with respect to a class of $L$-structures that has amalgamation and satisfies a few other simple conditions; see also
Ehrlich \cite{E1}. We would in any case still need Lemmas 3 and 4 to verify those
conditions.


\begin{thebibliography}{99}

\bibitem{AD}

M. Aschenbrenner  and L. van den Dries, \emph{H-fields and their Liouville extensions}, Math. Zeit. 242 (2002), pp. 543-548.

\bibitem{ADH1}

M. Aschenbrenner, L. van den Dries and J. van der Hoeven, \emph{Asymptotic Differential Algebra and Model Theory of Transseries}, Annals of Mathematics Studies, 195, Princeton University Press, Princeton, NJ (2017).

\bibitem{ADH2}\bysame,  \emph{Surreal numbers as a universal $H$-field}, J. Eur. Math. Soc. (forthcoming).

\bibitem{ADH3} \bysame, \emph{On Numbers, Germs, and Transseries}, to appear in the Proceedings of ICM 2018.

\bibitem{B}
A. Berarducci and V. Mantova, \emph{Surreal numbers, derivations and transseries}, J. Eur. Math. Soc. 20 (2018), pp. 339-390.

\bibitem{B1}\bysame, \emph{Transseries as germs of surreal functions}, Trans. Amer. Math. Soc., to appear, {\tt arXiv:1703.01995}.

\bibitem{C}
J. H. Conway, \emph{On Numbers and Games}, Academic Press, London (1976); \emph{Second Edition}, A K Peters, Ltd., Natick, Massachusetts (2001).

\bibitem{E1}
P. Ehrlich, \emph{Absolutely Saturated Models}, Fund. Math. 133 (1989), pp. 39-46.


\bibitem{J1}
B. J\'onsson, \emph{Universal relational systems}, Math. Scand. 4 (1956), pp. 193-208. 

\bibitem{J2}
\bysame,  \emph{Homogeneous universal relational systems}, Math. Scand. 8 (1960), pp. 137-142.


\bibitem{S}
M. C. Schmeling, \emph{Corps de transs\'eries}, Universit\'e Paris-VII, 2001.

\end{thebibliography}
\end{document}